\documentclass{amsart}

\usepackage{pstricks,pst-node}
\usepackage{amssymb,amsmath,amsthm}
\usepackage{hyperref}
\usepackage{doi,url}

\newtheorem{thm}{Theorem}[section]

\newtheorem{lem}{Lemma}[section]
\makeatletter
\renewcommand*{\c@lem}{\c@thm}
\makeatother

\newtheorem{prop}{Proposition}[section]
\makeatletter
\renewcommand*{\c@prop}{\c@thm}
\makeatother

\newtheorem{cor}{Corollary}[section]
\makeatletter
\renewcommand*{\c@cor}{\c@thm}
\makeatother

\theoremstyle{definition}

\newtheorem{example}{Example}[section]
\makeatletter
\renewcommand*{\c@example}{\c@thm}
\makeatother

\makeatletter
\renewcommand*{\c@conj}{\c@thm}
\makeatother

\theoremstyle{remark}
\newtheorem*{remark}{Remark}
\newtheorem*{note}{Note}

\newcount\show
\def\Label#1{\label{#1} \ifnum\show=1{\bf Label : [#1]}\fi}
\def\flr#1{\left\lfloor \frac{#1}{2}\right\rfloor}
\def\pp#1{\psi_{#1}}
\def\mat#1#2{\left(\begin{array}{#1}#2\end{array}\right)}
\def\seq#1#2#3{#1_1 #3 #1_2 #3 \def\temp{#3} \def\comma{,}
  \if\comma\temp \ldots \else \cdots \fi #3 #1_{#2} }
\def\d{\mathcal{D}}
\def\t{\mathcal{T}}
\def\Z{\mathbb{Z}}
\def\S{\mathfrak{S}}
\def\M{\mathfrak{M}_{n,m}^j}
\def\m{\mathfrak{m}_{n,m}^j}
\def\s{\mbox{\rm sign}}
\def\ao{\sigma_1}
\def\at{\sigma_2}
\def\bao{\bar{\sigma}_1}
\def\bat{\bar{\sigma}_2}
\def\lm{\lambda/\mu}
\def\ml{\mu/\lambda}
\def\ln{\lambda/\nu}
\def\mn{\mu/\nu}
\def\nm{\nu/\mu}
\def\dk{\delta_k}
\def\dr{\delta_r}
\def\tm{\tilde{\mu}}
\def\tl{\tilde{\lambda}}
\def\inverse{^{-1}}
\def\Inv{{\rm Inv}}
\def\inv{{\rm inv}}
\def\la{\lambda/\alpha}
\def\lb{\lambda/\beta}
\def\am{\alpha/\mu}
\def\bm{\beta/\mu}
\def\al{\alpha/\lambda}

\def\ta{\tilde{\alpha}}

\def\top{_{\rm top}}
\def\bot{_{\rm bottom}}
\def\lef{_{\rm left}}
\def\rig{_{\rm right}}
\def\st{^{\rm SDT}}
\def\up{\partial^+}
\def\down{\partial^-}
\def\sym{_{\rm sym}}

\def\qbinom#1#2{\genfrac{[}{]}{0pt}{}{#1}{#2}_q}
\def\qfactorial#1{\left[ #1 \right]_q!}

\def\wt#1{{\rm wt}_{#1}}
\def\ws{W} 
\def\G{G} 
\def\gdi{\Phi}

\def\led{\leq_{\rm d}}
\def\ld{<_{\rm d}}

\def\iff.{if and only if}
\def\dt.{domino tableau}
\def\dtl.{domino-tileable}
\def\skt.{skew tableau}
\def\SYT.{SYT}
\def\cp.{colored permutation}
\def\ci.{colored involution}
\def\RS.{Robinson-Schensted}
\def\si.{sign-imbalance}
\def\alg.{algorithm}
\def\iv.{involution}
\def\per.{permutation}
\def\sj.{Sj{\"o}strand}
\def\ss.{skew shape}
\def\gf.{generating function}
\def\gd.{growth diagram}
\def\rs.{reversed shape}

\newdimen\psunit
\newdimen\psraisebox
\def\psraise(#1,#2)#3{\psraisebox=\psunit \multiply\psraisebox -#1
  \divide\psraisebox #2 \divide\psraisebox 2 \advance\psraisebox 3.2pt
\raisebox{\psraisebox}{#3}}
\psset{unit=\psunit,linewidth=.4pt}
\newcount\ax \newcount\ay
\newcount\bx \newcount\by
\newcount\cx \newcount\cy
\newcount\dx \newcount\dy
\def\cell(#1,#2)[#3]{
\ax=#2 \ay=#1
\multiply\ay by-1
\bx=\ax \by=\ay
\cx=\ax \cy=\ay
\dx=\ax \dy=\ay
\advance\bx by-1
\advance\dy by1
\advance\cx by-1
\advance\cy by1
\pspolygon(\ax,\ay)(\bx,\by)(\cx,\cy)(\dx,\dy)(\ax,\ay)
\rput(\number\cx.5,
\ifnum\cy=0 -0.5\else\number\cy.5\fi){#3}
}

\begin{document}

\title{Skew domino Schensted correspondence and sign-imbalance}
\author{Jang Soo Kim} 
\thanks{The author was supported by the SRC program of Korea
  Science and Engineering Foundation (KOSEF) grant funded by the Korea
  government (MEST) (No. R11-2007-035-01002-0).}
\email{jskim@kaist.ac.kr}

\begin{abstract}
Using growth diagrams, we define a skew domino Schensted
correspondence which is a domino analogue of the skew
Robinson-Schensted correspondence due to Sagan and Stanley.  The
color-to-spin property of Shimozono and White is extended.  As an
application, we give a simple generating function for the weighted sum
of skew domino tableaux, which is a generalization of Stanley's
sign-imbalance formula.  The generating function gives a method to
calculate the generalized sign-imbalance formula. We also extend
Sj{\"o}strand's theorems on sign-imbalance of skew shapes.
\end{abstract}

\maketitle

\section{Introduction}

The domino Schensted correspondence is a bijection between \cp.s and
pairs of \dt.x of the same shape. It was first developed by Barbasch
and Vogan \cite{Barbasch1982} in 1982.  Garfinkle \cite{Garfinkle1992}
described this correspondence in terms of insertion.  Van Leeuwen
\cite{Leeuwen1996} described this correspondence using \gd.s and
extended it in the presence of a nonempty core. Shimozono and White
\cite{Shimozono2001} proved that this correspondence has the
color-to-spin property. Lam \cite{Lam2004} used \gd.s to prove the
color-to-spin property and identities involving \ci.s. Using these
properties, Lam \cite{Lam2004} obtained enumerative results for \dt.x
and proved Stanley's \si. conjectures \cite{Stanley2005}.

For a standard Young tableau (SYT) $T$, the sign of $T$ is defined to
be $\s(\pi)$, where $\pi$ is the \per. obtained by reading $T$ like a
book.  For example, if $T=$ \psraise (2,1){\pspicture (0,-2) (3,0)
  \cell(1,1)[1] \cell(1,2)[2] \cell(1,3)[4] \cell(2,1)[3]
  \cell(2,2)[5] \endpspicture} then $\s(T)=\s(12435)=-1$.  The
\si. $I_{\lambda}$ of a partition $\lambda$ is the sum of $\s(T)$ for
all SYTs $T$ of shape $\lambda$.  In \cite{Stanley2005}, Stanley
suggested the following interesting \si. formulas:
\begin{equation}\Label{one}
\sum_{\lambda \vdash
  n}x^{v(\lambda)}y^{h(\lambda)}z^{d(\lambda)}I_{\lambda}
=(x+y)^{{\flr n}},
\end{equation}
\begin{equation}\Label{two}
\sum_{\lambda \vdash n}(-1)^{v(\lambda)}I_{\lambda}^2 = 0,
\end{equation}
where $v(\lambda)$, $h(\lambda)$ and $d(\lambda)$ denote the maximum
numbers of vertical dominoes, horizontal dominoes and $2\times2$
rectangles respectively that can be placed in the Young diagram of
$\lambda$ without overlaps.

Reifegerste \cite{Reifegerste2004} and \sj. \cite{Sjostrand2005}
independently proved that if $\pi$ corresponds to $(P,Q)$ in the
\RS. correspondence and $sh(P)=\lambda$ then
\begin{equation}\Label{sjrei}
\s(\pi)=(-1)^{v(\lambda)} \s(P)\s(Q).
\end{equation}
Using Eq.~\eqref{sjrei}, Reifegerste \cite{Reifegerste2004} and \sj.
\cite{Sjostrand2005} proved Eq.~\eqref{two}.  \sj.
\cite{Sjostrand2005} also proved Eq.~\eqref{one} using Chess tableaux.

White \cite{White2001} observed that \si. is related to \dt.x and
proved that for a \dt. $D$,
\begin{equation}\Label{white}
\s(D)=(-1)^{ev(D)},  
\end{equation}
where $ev(D)$ is the number of vertical dominoes of $D$ in even
columns. 

Lam \cite{Lam2004} proved Eq.~\eqref{one} and Eq.~\eqref{two} using
\gd.s and Eq.~\eqref{white}.

There are some results about \si. for skew shapes.
\sj. \cite{Sjostrand2006} generalized Eq.~\eqref{sjrei} as follows: If
$sh(T)=sh(U)=\am$, $sh(P)=sh(Q)=\la$ and $\pi$ is an $n$-partial
permutation satisfying $(\pi,T,U)\leftrightarrow(P,Q)$ in the skew
\RS. correspondence of Sagan and Stanley \cite{Sagan1990}, then
\begin{equation}\Label{sjeq}
(-1)^{v(\lambda)}\s(P)\s(Q)=(-1)^{|\alpha|}(-1)^{v(\mu)+|\mu|}\s(T)\s(U)
\s(\overline{\pi}),  
\end{equation}
where $\overline{\pi}$ is the $n$-permutation extended from $\pi$ with
the smallest number of inversions. Using Eq.~\eqref{sjeq},
\sj. \cite{Sjostrand2006} generalized Eq.~\eqref{two} as follows: If
$\alpha$ is a fixed partition then
\begin{equation}\Label{sjsi}
\sum_{\la\vdash n}\!\!(-1)^{v(\lambda)} I_{\la}^2=
(-1)^n\!\!\!\sum_{\am\vdash n}\!\! (-1)^{v(\mu)} I_{\am}^2
+\frac{1-(-1)^n}2 \!\!\!\!\sum_{\am\vdash n-1}
\!\!\!\!\! (-1)^{v(\mu)} I_{\am}^2.
\end{equation}
Lam proved Eq.~\eqref{sjsi} once using signed differential posets
\cite{Lam2006a} and once, when $n$ is even, using the skew domino
Cauchy identity \cite{Lam2006}.

In this paper, inspired by Lam's work \cite{Lam2004}, we describe a
skew domino Schensted correspondence using \gd.s, which is a domino
analogue of the skew \RS. correspondence. This \gd. approach was used
in Roby's thesis \cite{Roby1991} to describe the skew
\RS. correspondence.  Fomin \cite{Fomin1995a} proved the existence of
the skew \RS. correspondence in a more general context using operators
on partitions.  The color-to-spin property and Lam's identities for
\ci.s are extended.  As an application, we generalize Eq.~\eqref{one}
to skew shapes.  We also generalize Eq.~\eqref{sjeq} to skew tableaux
$P$ and $Q$ of shape $\la$ and $\lb$ respectively, and then generalize
Eq.~\eqref{sjsi}.

We should note that, in the literature, there are two different
definitions of $\s(T)$ for a SYT $T$ of shape $\lm$. In
\cite{Lam2006,Sjostrand2006}, the sign of a SYT $T$ of shape $\lm$
does not consider the cells in $\mu$, but in \cite{Lam2006a}, it
does. However, if $sh(T)=sh(U)$ then the product $\s(T)\s(U)$ is the
same in both definitions, and so are Eq.~\eqref{sjeq} and
Eq.~\eqref{sjsi}.  In this paper, we use the definition of $\s(T)$ in
\cite{Lam2006a} and prove that Eq.~\eqref{white} still holds for skew
domino tableaux.

The rest of this paper is organized as follows.  In \autoref{pre}, we
define skew shapes, \rs.s, \dt.x and \cp.s.  In \autoref{growth}, we
introduce growth diagrams and a skew domino Schensted correspondence
and extend the color-to-spin property and Lam's identities for
\ci.s. We also find a \gf. for the weighted sum of \dt.x which turns
out to be closely related to \si..  In \autoref{si}, we define the
sign of a \skt.  and generalize Eq.~\eqref{one} to \ss.s.  The last
part of this section is devoted to finding a closed formula for
$\sum_{\lambda/\dk \vdash n} x^{v(\lambda/\dk)} y^{h(\lambda/\dk)}
z^{d(\lambda/\dk)} I_{\lambda/\dk}$, where $\dk=(k,k-1,\ldots,1)$. In
\autoref{sjthm}, we generalize Eq.~\eqref{sjeq} and Eq.~\eqref{sjsi}.

\section{Preliminaries}\Label{pre}
\subsection{Skew shapes and domino tableaux}
For a positive integer $n$, we denote $[n]=\{1,2,\ldots,n\}$. A {\em
  partition} $\lambda = (\seq{\lambda}{\ell}{,})$ of $n$, denoted by
$\lambda\vdash n$, is a weakly decreasing (possibly empty) sequence of
positive integers $\seq{\lambda}{\ell}{\geq}$ summing to $n$. Each
$\lambda_i$ is called the $i$-{\em th part} of $\lambda$. Let
$\ell(\lambda)$ denote the number of parts in $\lambda$.

A {\em cell} is a pair of positive integers.  The {\em Young diagram}
$Y(\lambda)$ of a partition $\lambda$ is the set of cells $(i,j)$ with
$i\leq \ell(\lambda)$ and $j\leq \lambda_i$.  We can draw the Young
diagram $Y(\lambda)$ by placing a square in the $i$-{\em th row} and
$j$-{\em th column} for each cell $(i,j)\in Y(\lambda)$. For example,
the drawing of the Young diagram of $(4,3,1)$ is $ \psraise
(3,1){\pspicture (0,-3) (4,0) \cell(1,1)[] \cell(1,2)[] \cell(1,3)[]
  \cell(1,4)[] \cell(2,1)[] \cell(2,2)[] \cell(2,3)[] \cell(3,1)[]
  \endpspicture} $. We will identify a partition $\lambda$ with its
Young diagram $Y(\lambda)$.

A {\em \ss.} $\lm$ is an ordered pair $(\lambda,\mu)$ of partitions
satisfying $\mu\subset\lambda$.  The {\em size} of $\lm$, denoted by
$|\lm|$, is the number of cells in $\lm$.  The notation $\lm\vdash n$
means that the size of $\lm$ is $n$. For example,
$(4,3,1)/(2,1)=
\psraise (3,1){\pspicture (0,-3) (4,0) \cell(1,3)[] \cell(1,4)[] \cell(2,2)[] \cell(2,3)[] \cell(3,1)[]\pspolygon[linestyle=dotted] (2,0) (2,-1) (1,-1) (1,-2) (0,-2) (0,0) \endpspicture}
$ is a \ss. of size $5$. 

A {\em domino} is a horizontal domino or a vertical domino where a
{\em horizontal} (resp. {\em vertical}) {\em domino} is a set of two
adjacent cells $(i,j)$ and $(i,j+1)$ (resp. $(i,j)$ and $(i+1,j)$).

A {\em (skew) standard Young tableau (\SYT.)} of shape $\lm\vdash n$
is a bijection $T$ from the set of cells in $\lm$ to $[n]$ such that
$T((i,j))\leq T((i',j'))$ whenever $i\leq i'$ and $j\leq j'$. For a
cell $c\in\lm$, we call the integer $T(c)$ the {\em entry} of $c$.  A
{\em (skew) standard domino tableau (SDT)} of shape $\lm\vdash 2n$ is
a \SYT. such that two cells with entries $2i-1$ and $2i$ make a domino
for each $i=1,2,\ldots,n$. Thus we can consider a SDT as a collection
of labeled dominoes.  For example, 
\psraise (3,1){\pspicture (0,-3) (4,0) \cell(1,2)[1] \cell(1,3)[2] \cell(1,4)[9] \cell(2,1)[3] \cell(2,2)[4] \cell(2,3)[7] \cell(2,4)[10] \cell(3,1)[5] \cell(3,2)[6] \cell(3,3)[8]\pspolygon[linestyle=dotted] (1,0) (1,-1) (0,-1) (0,0) \endpspicture}
 and
\psraise (3,1){\pspicture (0,-3) (4,0)\pspolygon[linestyle=dotted] (1,0) (1,-1) (0,-1) (0,0)\pspolygon (1,-1) (3,-1) (3,0) (1,0)\rput (2,-.5){1}\pspolygon (2,-1) (2,-2) (0,-2) (0,-1)\rput (1,-1.5){2}\pspolygon (2,-2) (2,-3) (0,-3) (0,-2)\rput (1,-2.5){3}\pspolygon (3,-1) (3,-3) (2,-3) (2,-1)\rput (2.5,-2){4}\pspolygon (3,-2) (4,-2) (4,0) (3,0)\rput (3.5,-1){5} \endpspicture}
 represent the same SDT. If there is a SDT of shape
$\lm$, then we say $\lm$ is {\em \dtl.}.

Let $\t(\lm)$ (resp. $\d(\lm)$) denote the set of all \SYT.s (resp. SDTs)
of shape $\lm$. Let $f^{\lm} = |\t(\lm)|$ and $d^{\lm}=|\d(\lm)|$.

For a given partition $\lambda$, let us take a maximal chain of
partitions $\lambda^{(m)}\subset \lambda^{(m-1)}\subset \cdots
\subset\lambda^{(0)}=\lambda$ such that
$\lambda^{(i-1)}/\lambda^{(i)}$ is a domino for $i=1,2,\ldots,
m$. Then the partition $\lambda^{(m)}$ is always the same and is
called the {\em $2$-core} of $\lambda$. We denote the $2$-core of
$\lambda$ by $\tl$.  Since there is no partition $\mu$ such that
$\tl/\mu$ is a domino, $\tl$ must be a {\em staircase partition}
$\dr=(r,r-1,\ldots,1)$ for some $r$. We refer the reader to
\cite{Garvan1990,James1981,Macdonald1995} for details of $p$-cores.

Let $v(\lm)$ (resp. $h(\lm)$) denote the number of cells in even
rows (resp. columns). Let $d(\lm)$ denote the number of cells both in
even columns and even rows.  It is easy to see that $h(\lambda)$,
$v(\lambda)$ and $d(\lambda)$ are the maximum numbers of horizontal
dominoes, vertical dominoes and $2\times 2$ rectangles respectively
that can be placed in $\lambda$ without overlaps.

For a SDT $D$, let $oh(D)$, $eh(D)$, $ov(D)$ and $ev(D)$ denote the
numbers of horizontal dominoes in odd rows, horizontal dominoes in
even rows, vertical dominoes in odd columns and vertical dominoes
in even columns respectively. The {\em spin} of a SDT is defined to
be the number of vertical dominoes divided by 2, that is,
$sp(D)=\frac12(ov(D)+ev(D))$.

Next we prove some relations between statistics of SDTs. We note that
these can also be proved by modifying Lam's results \cite{Lam2004} in
Proposition 14.

\begin{lem}\Label{stat}
If $D\in \d(\lm)$ then the following hold.
\begin{enumerate}
\item $oh(D)-eh(D) = \frac12 |\lm| - v(\lm)$
\item $ov(D)-ev(D) = \frac12 |\lm| - h(\lm)$
\item $eh(D)+ev(D) = d(\lm)$
\item $v(\lm) + h(\lm) = \frac12 |\lm| + 2\cdot d(\lm)$
\end{enumerate}
\end{lem}
\begin{proof}
Assign $1$ to the cells in odd rows and $-1$ to the cells in even rows
in $\lm$. Then the sum of all assigned numbers is $|\lm|-2v(\lm)$.
Each vertical domino contains both $1$ and $-1$.  Each horizontal
domino contains two $1$'s or two $-1$'s in accordance with the parity
of its row number.  Thus the sum is equal to $2oh(D)-2eh(D)$, which
proves the first identity. Similarly we can prove the second identity.

The right hand side of the third equation is the number of cells both
in even rows and even columns of $\lm$. A domino {\bf d} contains one
of these cells \iff. {\bf d} is either a horizontal domino in an even
row or a vertical domino in an even column.  Thus
$d(\lm)=eh(D)+ev(D)$.  

By the first three identities, we get the fourth:
\begin{align*}
  v(\lm)+h(\lm) &= |\lm| -(oh(D)+ov(D)-eh(D)-ev(D))\\
&= |\lm|-\left(\frac 12 |\lm| -2(eh(D)+ev(D))\right)\\
&= \frac 12 |\lm| + 2\cdot d(\lm). \qedhere
\end{align*}
\end{proof}

\begin{remark}
In \autoref{stat}, $(4)$ is not true if $\lm$ is not \dtl..
For example, if $\mu=(1)$ and $\lambda=(2,1)$, then $(4)$ does not
hold.
\end{remark}

\subsection{Reversed shapes}
Recall that a \ss. $\ml$ is a pair $(\mu,\lambda)$ of partitions with
$\lambda\subset \mu$. We define a {\em \rs.} $\lm$ to be a pair
$(\lambda,\mu)$ of partitions with $\lambda\subset\mu$ and denote
$\lm\vdash |\lambda|-|\mu|$.  Thus $\lm$ is a \rs. \iff. $\ml$ is a
\ss.. We also see that $\lm\vdash -n$ is equivalent to $\ml\vdash n$.
We extend each statistic $stat$ of \ss.s to \rs.s by defining
$stat(\lm)=-stat(\ml)$, i.e., $|\lm|=-|\ml|$, $v(\lm)=-v(\ml)$ and so
on. As a shape of a tableau, we will treat $\ml$ and $\lm$ equally,
that is, $\t(\lm)=\t(\ml)$ and $\d(\lm)=\d(\ml)$.

To avoid confusion we will always write a reversed shape with the
negative sign, that is, if we write $\lm\vdash n$ (resp.  $\lm\vdash
-n$) then it is always assumed that $n\geq0$ and $\lm$ is a skew
shape (resp. reversed shape).

The notion of \rs.s is not essential. However it will give us a simple
description for the generalization of Eq.~\eqref{one}. In
\autoref{si}, we will define the \si. $I_{\lm}$ of a \rs. $\lm\vdash
-2n$.

\subsection{Colored permutations and colored involutions}
A {\em \cp.} $\pi$ of $[n]$ is a permutation of $[n]$ equipped with an
assignment of bars to some integers. Let $\pi$ be a \cp..  The {\em
  total color} $tc(\pi)$ of $\pi$ is the number of barred
integers. The {\em permutation matrix} of $\pi$ is the matrix $M$ such
that $M(i,j)$ is equal to 1 if $\pi_i = j$; $-1$ if $\pi_i = \bar j$
and 0 otherwise.  Let $CP_n$ denote the set of \cp.s of $[n]$.

A \cp. $\pi$ is called an {\it \iv.} if the permutation matrix of
$\pi$ is symmetric. We denote the set of involutions in $CP_n$ by
$CI_n$. We will consider the empty word as an involution, thus
$CI_0=\{\emptyset\}$.  We can represent a \cp.  in cycle notation as
follows. Given a \cp. $\pi$, write the underlying permutation of $\pi$
in cycle notation, and put a bar over $i$ \iff. $i$ is barred in
$\pi$. For example, if $\pi=\bar341\bar52$ then
$\pi=(1\bar3)(24\bar5)$ in cycle notation.

Let $\pi$ be a \ci.. Then $\pi$ has only 1-cycles and 2-cycles, and
moreover, the two integers in a 2-cycle of $\pi$ are both barred or
both unbarred. Let $\ao(\pi)$, $\at(\pi)$, $\bao(\pi)$ and $\bat(\pi)$
denote the numbers of unbarred 1-cycles, unbarred 2-cycles, barred
1-cycles and barred 2-cycles in $\pi$ respectively.  For example, if
$\pi=(14)(\bar2)(\bar3\bar6)(5)(7)$ then $\ao(\pi)=2$, $\bao(\pi)=1$,
$\at(\pi)=1$ and $\bat(\pi)=1$.

We define the weight of a \ci. $\pi$ by
$$\wt{\pi}=\wt{\pi}(x,y,q)
=x^{\ao(\pi)}y^{\bao(\pi)}q^{\frac12{tc(\pi)}}.$$

Since a \ci. $\pi$ can be considered as a partition of $[n]$ into
1-subsets and 2-subsets with a possible bar on each subset and
$\frac12tc(\pi)= \frac12\bao(\pi)+\bat(\pi)$, by the exponential
formula \cite{Stanley1999}, we get the following exponential \gf.:
\begin{equation}\Label{wsi}
\sum_{n\geq0} \left(\sum_{\pi\in CI_n}\wt{\pi}\right)
  \frac{t^n}{n!}  =\exp\left((x + y\sqrt q)t+ (1+q)
  \frac{t^2}2\right).
\end{equation}

\section{Skew Domino Schensted Correspondence}\Label{growth}
\subsection{Definition of a \gd.}
In this section we introduce \gd.s.  Our definition is based on Lam's
 \cite{Lam2004}. We can define growth diagrams of an arbitrary skew
shape. Nevertheless, we will restrict our definition to rectangular
shapes for simplicity since we only need that case.  The reader is
referred to  \cite{Fomin1986, Fomin1994, Fomin1995} for details of
\gd.s.

For partitions $\lambda$ and $\mu$, we write $\mu\ld\lambda$ if $\lm$
is a domino and $\mu\led\lambda$ if $\mu=\lambda$ or $\mu\ld\lambda$.
A {\em d-chain} is a chain of partitions $\lambda^{(0)}\ld
\lambda^{(1)}\ld \cdots \ld \lambda^{(m)}$ and a {\em d-multichain} is
a multichain of partitions $\lambda^{(0)}\led \lambda^{(1)}\led \cdots
\led \lambda^{(m)}$. 

An $n\times m$ {\em growth array} $\Gamma$ is an array of partitions
$\Gamma_{(i,j)}$ for $0\leq i\leq n$ and $0\leq j\leq m$ such that
any two adjacent partitions are equal or differ by a domino, i.e.,
$\Gamma_{(i-1,j)}\led\Gamma_{(i,j)}$ and
$\Gamma_{(i,j-1)}\led\Gamma_{(i,j)}$ for all $1\leq i\leq n$ and
$1\leq j\leq m$.

An $n\times m$ {\em partial permutation matrix (PPM)} $M$ is an
$n\times m$ matrix whose elements are $1$, $-1$, or $0$, and which
contains at most one nonzero element in each row and column.  For a
PPM $M$, let $cp(M)$ denote the \cp. $\pi$ whose permutation matrix is
the matrix obtained from $M$ by removing the rows and columns
consisting of zeroes only.

An $n\times m$ {\em \gd.} $\G$ is a pair $(\Gamma, M)$, where
$\Gamma=\Gamma(G)$ is an $n\times m$ growth array and $M=M(G)$ is an
$n\times m$ PPM satisfying the following {\em local rules}.

Let $\nu = \Gamma_{(i-1,j-1)}$, $\mu=\Gamma_{(i-1,j)}$, $\rho =
\Gamma_{(i,j-1)}$ and $\lambda =\Gamma_{(i,j)}$. Then it must fall
into one of the following conditions which determine $\lambda$:
\begin{enumerate}
\item If $M(i,j)=1$ then $\nu=\mu=\rho$ and $\lambda$ is the
partition obtained from $\mu$ by adding a horizontal domino to
the first row.
\item If $M(i,j)=-1$ then $\nu=\mu=\rho$ and $\lambda$ is the partition
  obtained from $\mu$ by adding a vertical domino to the first column.
\item If $M(i,j)=0$ then there are five cases.
  \begin{enumerate}
  \item If $\nu=\mu$ or $\nu=\rho$ then $\lambda$ is the maximal
    partition among $\nu,\mu$ and $\rho$.
  \item If $\nu\ld\mu$, $\nu\ld\rho$, $\mu\ne\rho$ and $\mu/\nu
    \cap \rho/\nu =\emptyset$ then $\lambda=\mu\cup\rho$.
  \item If $\nu\ld\mu$, $\nu\ld\rho$, $\mu\ne\rho$ and $\mu/\nu
    \cap \rho/\nu \ne \emptyset$ then $\mn$ and $\rho/\nu$ share
    only one cell, say $(p,q)$, and $\lambda$ is the partition
    obtained from $\mu\cup\rho$ by adding the cell $(p+1,q+1)$.
  \item If $\nu\ld\mu$, $\nu\ld\rho$, $\mu=\rho$ and
    $\mn$ is a horizontal domino in $k$-th row then $\lambda$ is the
    partition obtained from $\mu$ by adding a horizontal domino to the
    $(k+1)$-th row.
  \item If $\nu\ld\mu$, $\nu\ld\rho$, $\mu=\rho$ and
    $\mn$ is a vertical domino in $k$-th column then $\lambda$ is the
    partition obtained from $\mu$ by adding a vertical domino to the
    $(k+1)$-th column.
  \end{enumerate}
\end{enumerate}

\begin{figure}
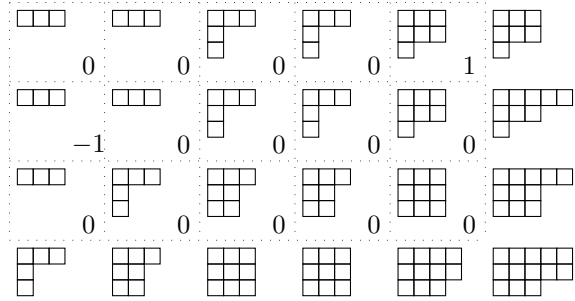

  \begin{center}
\psset{unit=6pt}
\psraise (18,1){\pspicture (0,-18) (36,0)\cell(1,1)[] \cell(1,2)[] \cell(1,3)[] \cell(1,7)[] \cell(1,8)[] \cell(1,9)[] \cell(1,13)[] \cell(1,14)[] \cell(1,15)[] \cell(2,13)[] \cell(3,13)[] \cell(1,19)[] \cell(1,20)[] \cell(1,21)[] \cell(2,19)[] \cell(3,19)[] \cell(1,25)[] \cell(1,26)[] \cell(1,27)[] \cell(2,25)[] \cell(2,26)[] \cell(2,27)[] \cell(3,25)[] \cell(1,31)[] \cell(1,32)[] \cell(1,33)[] \cell(2,31)[] \cell(2,32)[] \cell(2,33)[] \cell(3,31)[] \cell(6,1)[] \cell(6,2)[] \cell(6,3)[] \cell(6,7)[] \cell(6,8)[] \cell(6,9)[] \cell(6,13)[] \cell(6,14)[] \cell(6,15)[] \cell(7,13)[] \cell(8,13)[] \cell(6,19)[] \cell(6,20)[] \cell(6,21)[] \cell(7,19)[] \cell(8,19)[] \cell(6,25)[] \cell(6,26)[] \cell(6,27)[] \cell(7,25)[] \cell(7,26)[] \cell(7,27)[] \cell(8,25)[] \cell(6,31)[] \cell(6,32)[] \cell(6,33)[] \cell(6,34)[] \cell(6,35)[] \cell(7,31)[] \cell(7,32)[] \cell(7,33)[] \cell(8,31)[] \cell(11,1)[] \cell(11,2)[] \cell(11,3)[] \cell(11,7)[] \cell(11,8)[] \cell(11,9)[] \cell(12,7)[] \cell(13,7)[] \cell(11,13)[] \cell(11,14)[] \cell(11,15)[] \cell(12,13)[] \cell(12,14)[] \cell(13,13)[] \cell(13,14)[] \cell(11,19)[] \cell(11,20)[] \cell(11,21)[] \cell(12,19)[] \cell(12,20)[] \cell(13,19)[] \cell(13,20)[] \cell(11,25)[] \cell(11,26)[] \cell(11,27)[] \cell(12,25)[] \cell(12,26)[] \cell(12,27)[] \cell(13,25)[] \cell(13,26)[] \cell(13,27)[] \cell(11,31)[] \cell(11,32)[] \cell(11,33)[] \cell(11,34)[] \cell(11,35)[] \cell(12,31)[] \cell(12,32)[] \cell(12,33)[] \cell(13,31)[] \cell(13,32)[] \cell(13,33)[] \cell(16,1)[] \cell(16,2)[] \cell(16,3)[] \cell(17,1)[] \cell(18,1)[] \cell(16,7)[] \cell(16,8)[] \cell(16,9)[] \cell(17,7)[] \cell(17,8)[] \cell(18,7)[] \cell(18,8)[] \cell(16,13)[] \cell(16,14)[] \cell(16,15)[] \cell(17,13)[] \cell(17,14)[] \cell(17,15)[] \cell(18,13)[] \cell(18,14)[] \cell(18,15)[] \cell(16,19)[] \cell(16,20)[] \cell(16,21)[] \cell(17,19)[] \cell(17,20)[] \cell(17,21)[] \cell(18,19)[] \cell(18,20)[] \cell(18,21)[] \cell(16,25)[] \cell(16,26)[] \cell(16,27)[] \cell(16,28)[] \cell(17,25)[] \cell(17,26)[] \cell(17,27)[] \cell(17,28)[] \cell(18,25)[] \cell(18,26)[] \cell(18,27)[] \cell(16,31)[] \cell(16,32)[] \cell(16,33)[] \cell(16,34)[] \cell(16,35)[] \cell(17,31)[] \cell(17,32)[] \cell(17,33)[] \cell(17,34)[] \cell(17,35)[] \cell(18,31)[] \cell(18,32)[] \cell(18,33)[] \rput (4.5,-3.5){$0$}\rput (10.5,-3.5){$0$}\rput (16.5,-3.5){$0$}\rput (22.5,-3.5){$0$}\rput (28.5,-3.5){$1$}\rput (4.5,-8.5){$-1$}\rput (10.5,-8.5){$0$}\rput (16.5,-8.5){$0$}\rput (22.5,-8.5){$0$}\rput (28.5,-8.5){$0$}\rput (4.5,-13.5){$0$}\rput (10.5,-13.5){$0$}\rput (16.5,-13.5){$0$}\rput (22.5,-13.5){$0$}\rput (28.5,-13.5){$0$}\psline[linestyle=dotted] (-.5,.5) (29.5,.5)\psline[linestyle=dotted] (-.5,-4.5) (29.5,-4.5)\psline[linestyle=dotted] (-.5,-9.5) (29.5,-9.5)\psline[linestyle=dotted] (-.5,-14.5) (29.5,-14.5)\psline[linestyle=dotted] (-.5,.5) (-.5,-14.5)\psline[linestyle=dotted] (5.5,.5) (5.5,-14.5)\psline[linestyle=dotted] (11.5,.5) (11.5,-14.5)\psline[linestyle=dotted] (17.5,.5) (17.5,-14.5)\psline[linestyle=dotted] (23.5,.5) (23.5,-14.5)\psline[linestyle=dotted] (29.5,.5) (29.5,-14.5) \endpspicture}
\psset{unit=10pt}
  \end{center}
\caption{A $3\times5$ growth diagram.}\label{fig:growth}
\end{figure}

For example, see Fig.~\ref{fig:growth}, which represents a
\gd. $G=(\Gamma,M)$ with
$M=\mat{ccccc}{0&0&0&0&1\\-1&0&0&0&0\\0&0&0&0&0}$.

\subsection{Skew domino Schensted correspondence}\label{subsec}
Let $C=\lambda^{(0)}\led \lambda^{(1)}\led \cdots \led \lambda^{(n)}$
be a d-multichain. We can naturally construct a SDT from $C$ as
follows. Let $({\bf d}_1,{\bf d}_2,\ldots, {\bf d}_k)$ be the sequence
of dominoes obtained by removing the empty skew shapes (if any) from
$(\lambda^{(1)}/\lambda^{(0)}, \lambda^{(2)}/\lambda^{(1)},\ldots,
\lambda^{(n)}/\lambda^{(n-1)})$. Then $C\st$ denotes the SDT of shape
$\lambda^{(n)}/\lambda^{(0)}$ whose domino with entry $i$ is ${\bf
  d}_i$ for $i=1,2,\ldots,k$.

Let $G=(\Gamma,M)$ be an $n\times m$ \gd..  We define four special
d-multichains of $G$ as follows:
\begin{align*}
\G\top &= \Gamma_{(0,0)} \led \Gamma_{(0,1)} \led \cdots \led \Gamma_{(0,m)},\\
\G\bot &= \Gamma_{(n,0)} \led \Gamma_{(n,1)} \led \cdots \led \Gamma_{(n,m)},\\
\G\lef &= \Gamma_{(0, 0)} \led \Gamma_{(1,0)} \led \cdots \led \Gamma_{(n,0)},\\
\G\rig &= \Gamma_{(0,m)} \led \Gamma_{(1,m)} \led \cdots \led \Gamma_{(n,m)}. 
\end{align*}
If both $G\bot$ and $G\rig$ are d-chains, then we call $G$ a {\em full
  \gd.}.  The local rules say that $G$ is completely determined by
$G\top$, $G\lef$ and $M$. On the other hand, one can easily see that
the local rules are invertible in the sense that $\Gamma_{(i-1,j-1)}$
and $M(i,j)$ are determined by $\Gamma_{(i-1,j)}$, $\Gamma_{(i,j-1)}$
and $\Gamma_{(i,j)}$. Thus a \gd. $G$ is also completely determined by
$G\bot$ and $G\rig$. We define $\up(G)$ to be the pair
$(G\bot\st,\G\rig\st)$, and $\down(G)$ to be the triple
$(G\top\st,G\lef\st,M)$. For example, if $G$ is the the \gd. in
Fig.~\ref{fig:growth}, then
\begin{align*}
\up(G) &= \left(\:
\psraise (3,1){\pspicture (0,-3) (5,0)\pspolygon[linestyle=dotted] (3,0) (3,-1) (1,-1) (1,-3) (0,-3) (0,0)\pspolygon (2,-1) (2,-3) (1,-3) (1,-1)\rput (1.5,-2){1}\pspolygon (3,-1) (3,-3) (2,-3) (2,-1)\rput (2.5,-2){2}\pspolygon (3,-2) (4,-2) (4,0) (3,0)\rput (3.5,-1){3}\pspolygon (4,-2) (5,-2) (5,0) (4,0)\rput (4.5,-1){4} \endpspicture}
\:,
\psraise (3,1){\pspicture (0,-3) (5,0)\pspolygon[linestyle=dotted] (3,0) (3,-2) (1,-2) (1,-3) (0,-3) (0,0)\pspolygon (3,-1) (5,-1) (5,0) (3,0)\rput (4,-.5){1}\pspolygon (3,-2) (3,-3) (1,-3) (1,-2)\rput (2,-2.5){2}\pspolygon (3,-2) (5,-2) (5,-1) (3,-1)\rput (4,-1.5){3} \endpspicture}
\:\right), \\
\down(G) &= \left(\:
\psraise (3,1){\pspicture (0,-3) (3,0)\pspolygon[linestyle=dotted] (3,0) (3,-1) (0,-1) (0,0)\pspolygon (1,-1) (1,-3) (0,-3) (0,-1)\rput (.5,-2){1}\pspolygon (3,-1) (3,-2) (1,-2) (1,-1)\rput (2,-1.5){2} \endpspicture}
\:,
\psraise (3,1){\pspicture (0,-3) (3,0)\pspolygon[linestyle=dotted] (3,0) (3,-1) (0,-1) (0,0)\pspolygon (1,-1) (1,-3) (0,-3) (0,-1)\rput (.5,-2){1} \endpspicture}
\:,
\mat{ccccc}{ 0&0&0&0&1\\-1&0&0&0&0\\0&0&0&0&0} \:\right).
\end{align*}

Let $\alpha$ and $\beta$ be partitions. Let
$\mathfrak{G}_{n,m}^{\alpha,\beta}$ denote the set of all $n\times m$
full \gd.s $G=(\Gamma,M)$ satisfying $\Gamma_{(n,0)}=\alpha$ and
$\Gamma_{(0,m)}=\beta$. 

Let $\M$ denote the set of all $n\times m$ PPMs (partial permutation
matrices) with $j$ nonzero elements.

\begin{lem}\Label{full}
Let $C=C_{(0)}\led C_{(1)}\led \cdots \led C_{(m)}$ and
$C'=C'_{(0)}\led C'_{(1)}\led \cdots \led C'_{(n)}$ be d-multichains.
Let $M$ be an $n\times m$ PPM.  Then, there is a (necessarily unique)
$n\times m$ full \gd. $G=(\Gamma,M)$ such that $G\top=C$ and
$G\lef=C'$ \iff.  the following conditions hold:
\begin{enumerate}
\item For $1\leq i\leq n$, the $i$-th row of $M$ contains a
  nonzero element \iff. $C'_{(i-1)}=C'_{(i)}$.
\item For $1\leq j\leq m$, the $j$-th column of $M$ contains a
  nonzero element \iff. $C_{(j-1)}=C_{(j)}$.
\end{enumerate}
\end{lem}
\begin{proof}
We can check this easily by the local rules.
\end{proof}

\begin{lem}
Let $\alpha$ and $\beta$ be partitions. Then the maps $\up$ and
$\down$ induce the following bijections:
\begin{align*}
  \up &: \mathfrak{G}_{n,m}^{\alpha,\beta} \rightarrow
  \bigcup_{\substack{\la\vdash 2m \\ \lb\vdash 2n}}
  \d(\la)\times\d(\lb),\\ \down &: \mathfrak{G}_{n,m}^{\alpha,\beta}
  \rightarrow \bigcup_{j\geq0} \left( \bigcup_{\substack{\bm\vdash
      2(m-j) \\ \am\vdash 2(n-j)}} \d(\bm)\times\d(\am)\times
  \M \right).
\end{align*}
\end{lem}
\begin{proof}
By the local rules, a \gd. $G=(\Gamma,M)$ is determined by the pair
$(G\bot,G\rig)$ or the triple $(G\top, G\lef, M)$.  Moreover, if $G$
is a full \gd., then $(G\bot,G\rig)$ and $(G\top, G\lef, M)$ are in
bijection with $(G\bot\st,G\rig\st)$ and $(G\top\st, G\lef\st, M)$
respectively by \autoref{full}.  Thus $\up$ and $\down$ are invertible
for full \gd.s. The surjectiveness of $\up$ and $\down$ follows from
the local rules and \autoref{full}.
\end{proof}

Now we get a skew domino Schensted correspondence.
\begin{thm}\Label{gdi}
Let $\alpha$ and $\beta$ be fixed partitions and $n$ and $m$ be fixed
nonnegative integers.  Then $\gdi=\up\circ(\down)\inverse$
induces a bijection
$$\gdi : \bigcup_{j\geq0} \left( \bigcup_{\substack{\bm\vdash 2(m-j)
    \\ \am\vdash 2(n-j)}} \d(\bm)\times\d(\am)\times \M
\right) \rightarrow \bigcup_{\substack{\la\vdash 2m \\ \lb\vdash 2n}}
\d(\la)\times\d(\lb).$$
\end{thm}
We note that if $\pi$ corresponds to $(P,Q)$ in the domino Schensted
correspondence with the core $\dr$ and $M$ is the permutation matrix
of $\pi$ then $\gdi(\emptyset_{\dr}, \emptyset_{\dr}, M) = (P,Q)$,
where $\emptyset_{\dr}$ is the empty SDT of shape $\dr/\dr$.  The
bijection $\gdi$ is a domino analogue of the skew \RS. correspondence,
which was first developed using external and internal insertion by
Sagan and Stanley \cite{Sagan1990} and was interpreted in terms of
\gd.s, as we did here, by Roby \cite{Roby1991}. Fomin
\cite{Fomin1995a} proved the existence of the skew \RS. correspondence
using operators on partitions.

Since the local rules are symmetric we get the following proposition
immediately.
\begin{prop}
Let $\gdi(U,V,M)=(P,Q)$. Then $\gdi(V,U,M^T)=(Q,P)$.
\end{prop}
In the above proposition, if $U=V$ and $M$ is symmetric then
$\gdi(U,U,M)=(P,P)$. Let $\gdi\sym(U,M)=P$.  Then we get another
bijection.
\begin{cor}\Label{sgdi}
Let $\alpha$ be a fixed partition and $n$ be a fixed nonnegative
integer.  Then $\gdi\sym$ induces a bijection
$$\gdi\sym : \bigcup_{j\geq0} \left( \bigcup_{\am\vdash 2(n-j)}
\d(\am)\times \mathfrak{SM}_{n}^j \right) \rightarrow
\bigcup_{\la\vdash 2n} \d(\la),$$ where $\mathfrak{SM}_{n}^j$ denotes
the set of all symmetric $n\times n$ PPMs with $j$ nonzero elements.
\end{cor}

Shimozono and White \cite{Shimozono2001} proved that the domino
Schensted correspondence has the color-to-spin property, that is, if
$\pi$ corresponds to $(P,Q)$ then $tc(\pi)=sp(P)+sp(Q)$. The next
proposition generalizes this property. The proof is the same as Lam's
\cite[Lemma 8]{Lam2004}.
\begin{prop}\Label{colortospin}
Let $\gdi(U,V,M)=(P,Q)$ and $\pi=cp(M)$. Then
$$tc(\pi)=sp(P)+sp(Q)-sp(U)-sp(V).$$
\end{prop}
\begin{proof}
By the local rules, we can check that the following value is 1 if
$M(i,j)=-1$ and 0 otherwise: $sp(\Gamma_{(i,j)}/\Gamma_{(i-1,j)})+
sp(\Gamma_{(i,j)}/\Gamma_{(i,j-1)})-
sp(\Gamma_{(i-1,j)}/\Gamma_{(i-1,j-1)})-
sp(\Gamma_{(i,j-1)}/\Gamma_{(i-1,j-1)})$. By adding up these for all
$1\leq i\leq n$ and $1\leq j\leq m$, we finish the proof.
\end{proof}

Lam \cite{Lam2004} proved that if a \ci. $\pi$ corresponds to $(D,D)$
in the domino Schensted correspondence then $\bao(\pi)=ov(D)-ev(D)$
and $\bat(\pi)=ev(D)$. We can generalize Lam's results.
\begin{prop}\label{wt}
Let $M$ be an $n\times n$ symmetric PPM and $\pi=cp(M)$.  Let $U$ and
$D$ be SDTs satisfying $\gdi\sym(U,M)=D$. Then we have
  \begin{align*}
    \ao(\pi) &= (oh(D)-eh(D))+(oh(U)-eh(U)),\\
    \bao(\pi) &= (ov(D)-ev(D))+(ov(U)-ev(U)),\\
    \at(\pi) &=  eh(D) - oh(U),\\
    \bat(\pi) &=  ev(D) - ov(U).
  \end{align*}
\end{prop}
\begin{proof}
We will prove the second and the fourth identities. The remaining can
be proved similarly. By \autoref{colortospin}, we have
$$\bao(\pi) +2\bat(\pi) = (ov(D)+ev(D)) -
(ov(U)+ev(U)).$$ Thus, it is sufficient to show that
$$\bao(\pi) = (ov(D)-ev(D))+(ov(U)-ev(U)).$$ Let $G=(\Gamma,M)$ be the
corresponding $n\times n$ \gd. $(\up)\inverse (D,D)$ and let
$\nu=\Gamma_{(0,0)}$, $\mu=\Gamma_{(n,0)}=\Gamma_{(0,n)}$ and
$\lambda=\Gamma_{(n,n)}$. Then by \autoref{stat},
\begin{align*}
(ov(D)-ev(D))+(ov(U)-ev(U)) &= \frac{|\lambda/\mu|}2 - h(\lambda/\mu)
  + \frac{|\mn|}2 - h(\mn)\\ &=\frac{|\lambda/\nu|}2 - h(\lambda/\nu).
\end{align*}
One can check that $\frac12 |\Gamma_{(i,i)}/\Gamma_{(i-1,i-1)}| -
h(\Gamma_{(i,i)}/\Gamma_{(i-1,i-1)})$ is 1 if $M(i,i)=-1$ and 0
otherwise. By adding up these for all $1\leq i\leq n$, we finish
the proof.
\end{proof}

As an application of our skew domino Schensted correspondence, we get
some enumerative results.  Following Lam's notation \cite{Lam2004},
let
$$ f^{\lm}_2(q) = \sum_{D\in\d(\lm)} q^{sp(D)}.$$
\begin{cor}\Label{double}
Let $\alpha$ and $\beta$ be fixed partitions and $n$ and $m$ be fixed
nonnegative integers. Then,
$$\sum_{\substack{\la\vdash 2m \\ \lb\vdash 2n}} f^{\la}_2(q)
f^{\lb}_2(q) = \sum_{j\geq0} \binom{n}{j}\binom{m}{j} (1+q)^j j!
\sum_{ \substack{\bm\vdash 2(m-j) \\ \am \vdash 2(n-j)}} f^{\bm}_2(q)
f^{\am}_2(q).$$
\end{cor}
\begin{proof}
This is immediate from \autoref{gdi} and \autoref{colortospin}
and the following identity: $\sum_{\pi\in CP_j} q^{tc(\pi)}=(1+q)^j j!$.
\end{proof}
For a partition $\lambda=(\seq{\lambda}{l}{,})$, let $2\lambda$ denote
the partition $(\seq{2\lambda}{l}{,})$.  We can consider a \SYT. of
shape $\lambda/\mu$ as a SDT of shape $2\lambda/2\mu$ consisting of
horizontal dominoes by identifying a cell with a horizontal domino.

There are three interesting specializations of \autoref{double}.  When
$q=0$ in \autoref{double}, we get the following corollary due to Sagan
and Stanley \cite{Sagan1990}. We note that Roby \cite{Roby1991} also
proved the following corollary using \gd.s and our proof is
essentially the same as Roby's.
\begin{cor} \cite[Sagan and Stanley]{Sagan1990} \label{thm:ss}
Let $\alpha$ and $\beta$ be fixed partitions and $n$ and $m$ be fixed
nonnegative integers. Then,
$$\sum_{\substack{\la\vdash m \\ \lb\vdash n}} f^{\la} f^{\lb} =
\sum_{j\geq0} \binom{n}{j}\binom{m}{j} j!
\sum_{\substack{\bm\vdash m-j \\ \am \vdash n-j}} f^{\bm} f^{\am}.$$
\end{cor}

When we set $q=1$ in \autoref{double}, we get a domino analogue.
\begin{cor}
Let $\alpha$ and $\beta$ be fixed partitions and $n$ and $m$ be fixed
nonnegative integers. Then,
$$\sum_{\substack{\la\vdash 2m \\ {\lb\vdash 2n}}} d^{\la} d^{\lb} =
\sum_{j\geq0} \binom{n}{j}\binom{m}{j} 2^j j!
\sum_{\substack{\bm\vdash 2(m-j) \\ \am \vdash 2(n-j)}} d^{\bm}d^{\am}.$$
\end{cor}

If $q=-1$ then, as we will see in the next section, \autoref{double}
induces a sign-imbalance formula.

\subsection{The weighted sum of domino tableaux}
For a SDT $D$, we define the weight $\wt{D}$ of $D$ by
$$\wt{D}=\wt{D}(x,y,q) = x^{oh(D)-eh(D)}y^{ov(D)-ev(D)}q^{sp(D)}.$$
Note that if a \ci. $\pi$ corresponds to $(D,D)$ in the domino Schensted
correspondence, then $\wt{D}=\wt{\pi}$.

Recall that a \rs. $\lm\vdash -n$ is the one obtained by reversing a
\ss. $\ml\vdash n$. For a reversed shape $\lm\vdash -n$, we define
$f_2^{\lm}(q)=f_2^{\ml}(q)$.

For a fixed partition $\alpha$ and an integer $n\geq0$, we define
$$\ws_n^{\alpha} = \ws_n^{\alpha}(x,y,q) = \sum_{\la\vdash 2n} x^{\frac12|\la|-v(\la)}
y^{\frac12|\la|-h(\la)} f_2^{\la}(q),$$ 
and
$$\ws_{-n}^{\alpha} = \ws_{-n}^{\alpha}(x,y,q) = \sum_{\la\vdash -2n} 
x^{\frac12|\la|-v(\la)} y^{\frac12|\la|-h(\la)} f_2^{\la}(q).$$ Then,
$$\ws_{-n}^{\alpha} = \sum_{\al\vdash 2n} x^{-n+v(\al)} y^{-n+h(\al)}
f_2^{\al}(q).$$ Thus $W_{-n}^{\alpha}=0$ if $\alpha/\ta\vdash 2k$ and $n>k$.

By \autoref{stat}, if $n\geq0$ then $\ws_n^{\alpha}$ is
the weighted sum of certain SDTs:
$$\ws_n^{\alpha} = \sum_{\la\vdash 2n}\sum_{D\in\d(\la)} \wt{D}.$$ 

We note that $\ws_n^{\alpha}$ is a modified generalization of $h_r(n)$ in
Lam's paper~ \cite{Lam2004}:
$$h_r(n) = \sum_{\lambda/\dr\vdash 2n} a^{(o(\lambda)-o(\dr))/2}
b^{(o(\lambda')-o(\dr))/2} c^{d(\lambda)-d(\dr)} f_2^{\lambda}(q),$$
where $o(\lambda)$ denotes the number of odd parts in $\lambda$.  One
can check that $h_r(n)=c^{\frac n2} \ws_n^{\dr}( bc^{-\frac12},
ac^{-\frac12}, q)$.

\begin{thm}\Label{weightedsum}
  Let $\alpha$ be a fixed partition with $\alpha/\ta \vdash 2k$ and
  $n\geq0$. Then
$$\ws_n^{\alpha} = \sum_{j=0}^k \binom{n}{j} \ws_{-j}^{\alpha} \sum_{\pi\in
  CI_{n-j}} \wt{\pi}.$$
\end{thm}
\begin{proof}
Let $\la\vdash 2n$ and $D\in \d(\la)$.  Recall the bijection
$\gdi\sym$ in \autoref{sgdi}. Let $(\gdi\sym)\inverse(D)=(U,M)$,
$sh(U)=\am\vdash 2j$ and $cp(M)=\pi$. Then $j\leq k$ and $\pi\in
CI_{n-j}$. By \autoref{colortospin}, \autoref{wt} and \autoref{stat},
we have
$$\wt{D} = x^{-j+v(\am)} y^{-j+h(\am)} q^{sp(U)} \wt{\pi}.$$ Since $M$
is determined by $\pi$ and choosing $j$ nonzero rows, by
\autoref{sgdi},
\begin{align*}
\ws_n^{\alpha} &= \sum_{\la\vdash 2n} \sum_{D\in\d(\la)} \wt{D}\\ &=
\sum_{j=0}^k \binom{n}{j} \sum_{\am\vdash 2j} x^{-j+v(\am)}
y^{-j+h(\am)} f_2^{\am}(q) \sum_{\pi\in CI_{n-j}} \wt{\pi}\\ &=
\sum_{j=0}^k \binom{n}{j} \ws_{-j}^{\alpha} \sum_{\pi\in CI_{n-j}}
\wt{\pi}. \qedhere
\end{align*}
\end{proof}

Using \autoref{weightedsum} and Eq.~\eqref{wsi}, we get a simple
\gf. for the weighted sum.
\begin{cor}\Label{gf}
Let $\alpha$ be a fixed partition. Then
$$\frac{\sum_{n\geq0}\ws_n^{\alpha} t^n/n!}  {\sum_{n\geq0}
  \ws_{-n}^{\alpha} t^n/n!}  = \exp\left( (x+y\sqrt q)t +
(1+q)\frac{t^2}{2} \right).$$
\end{cor}
If we substitute $\alpha$, $x$, $y$ and $t$ in \autoref{gf} with
$\dr$, $bc^{-\frac12}$, $ac^{-\frac12}$ and $c^{\frac12} t$
respectively then we get Lam's result \cite{Lam2004}:
$$\sum_{n\geq0}h_r(n) \frac{t^n}{n!} = \exp\left( (b+a\sqrt q)t +
  c(1+q)\frac{t^2}{2} \right).$$

By the argument following \autoref{double}, if we set $x=1$ and
$y=q=0$ in \autoref{weightedsum}, then we obtain Sagan and Stanley's
theorem \cite{Sagan1990} which was reproved by Roby \cite{Roby1991},
Stanley \cite{Stanley2003} and Jaggard \cite{Jaggard2005}.

\begin{cor} \cite[Sagan and Stanley]{Sagan1990}\label{thm:ss1}
Let $\alpha\vdash k$ be a fixed partition. Then
$$\sum_{\la\vdash n}f^{\la} = \sum_{j=0}^{k} \binom{n}{j}t_{n-j}
\sum_{\am \vdash j} f^{\am} ,$$ where $t_m$ denotes the number of
\iv.s of $[m]$.
\end{cor}

If we set $x=y=q=1$ in \autoref{weightedsum}, we get the following
domino analogue.
\begin{cor}
Let $\alpha$ be a fixed partition with $\alpha /\ta\vdash 2k$. Then
$$\sum_{\la\vdash 2n}d^{\la} = \sum_{j=0}^{k} \binom{n}{j} \xi_{n-j}
\sum_{\am \vdash 2j} d^{\am},$$ where $\xi_m$ denotes the number of
\ci.s of $[m]$.
\end{cor}

\begin{note}
The skew Cauchy identities corresponding to \autoref{thm:ss} and
\autoref{thm:ss1} were obtained by Zelvinsky in the 1985 translation
of \cite{Macdonald1995}.  The skew domino Cauchy identity
corresponding to \autoref{double} was introduced by Lam
\cite{Lam2006}.
\end{note}

\section{A generalized sign-imbalance formula}\Label{si}
\subsection{Definition of the sign of a skew SYT}
For two cells $a=(i,j)$ and $b=(i',j')$, we write $a\lhd b$ if $i<i'$
or ($i=i'$ and $j<j'$).  For a \SYT. $T$, we denote $\Inv(T)=\{
(a,b):a\lhd b,T(a)>T(b)\}$ and $\inv(T)=|\Inv(T)|$.

The sign of a SYT $T$ is defined by $\s(T)=(-1)^{\inv (T)}$. The {\em
  \si.} $I_{\lambda}$ of a partition $\lambda$ is defined by
$$I_{\lambda}=\sum_{T\in\t(\lambda)} \s(T).$$ 

The purpose of this section is to define $I_{\lm}$ and generalize
Eq.~\eqref{one} and Eq.~\eqref{two}.  In the literature, there are two
different definitions for the sign of a skew \SYT. $T$. We will write
them as $\s_1(T)$ and $\s_2(T)$ temporarily.
\sj. \cite{Sjostrand2006} and Lam \cite{Lam2006} used $\s_1(T)$
defined by
$$\s_1(T)=(-1)^{\inv(T)}.$$ Lam \cite{Lam2006a} used $\s_2(T)$, which
we will use in this paper.  To define $\s_2(T)$, we introduce an
operation on two \SYT.s.

Assume $\mu\subset\nu\subset\lambda$ and $\nu/\mu\vdash k$. Let $T_1$
and $T_2$ be \SYT.s of shape $\nu/\mu$ and $\ln$ respectively. Then we
define $T_1\diamond T_2$ to be the \SYT. $T\in\t(\lm)$ such
that
$$T(c)=\left\{
\begin{array}{ll}
T_1(c),&\mbox{if $c\in \nu/\mu$},\\ T_2(c)+k,&\mbox{if $c\in \ln$}.
\end{array}\right.$$
For example, if $T_1=$ 
\psraise (3,1){\pspicture (0,-3) (3,0) \cell(1,3)[1] \cell(2,1)[2] \cell(2,2)[3]\pspolygon[linestyle=dotted] (2,0) (2,-1) (0,-1) (0,0) \endpspicture}
 and $T_2=$ 
\psraise (3,1){\pspicture (0,-3) (4,0) \cell(1,4)[4] \cell(2,3)[2] \cell(3,1)[1] \cell(3,2)[3]\pspolygon[linestyle=dotted] (3,0) (3,-1) (2,-1) (2,-2) (0,-2) (0,0) \endpspicture}
then $T_1\diamond T_2=$ 
\psraise (3,1){\pspicture (0,-3) (4,0) \cell(1,3)[1] \cell(1,4)[7] \cell(2,1)[2] \cell(2,2)[3] \cell(2,3)[5] \cell(3,1)[4] \cell(3,2)[6]\pspolygon[linestyle=dotted] (2,0) (2,-1) (0,-1) (0,0) \endpspicture}.

Now we define $\s_2(T)$ for a \SYT. $T\in\t(\lm)$ by
$$\s_2(T)=\s(T_0)\s(T_0\diamond T),$$ where $T_0$ is an arbitrary SYT
of shape $\mu$. It is straightforward to show the next proposition
which implies that $\s_2$ is well-defined.
\begin{prop}\Label{sign}
Let $T$ be a \SYT. of shape $\lm$. Then $\s(T_0)\s(T_0\diamond T)$ is
independent of the choice of $T_0\in \t(\mu)$.  Moreover,
$\s(T_0)\s(T_0\diamond T)=(-1)^m\s_1(T)$, where
$m=\sum_{i\geq1}(\lambda_i-\mu_i)\cdot \sum_{j>i}\mu_j$.
\end{prop}

We take $\s_2(T)$ for the sign of a \SYT. $T$. From now on, we will
write $\s(T)$ instead of $\s_2(T)$.  The $\s(T)$ has the following
product property.
\begin{prop}\Label{prod}
Let $\mu\subset\nu\subset\lambda$, $T_1\in\t(\nm)$ and
$T_2\in\t(\ln)$. Then 
$$\s(T_1\diamond T_2)=\s(T_1)\s(T_2).$$
\end{prop}
\begin{proof}
Let $T$ be a SYT of shape $\mu$. Then $T\diamond T_1$ is a SYT. Thus
\begin{align*}
  \s(T_1)\s(T_2) &= \s(T) \s(T\diamond T_1) \s(T\diamond T_1)
  \s(T\diamond T_1 \diamond T_2)\\ &= \s(T) \s(T\diamond T_1 \diamond
  T_2) = \s(T_1 \diamond T_2). \qedhere
  \end{align*}
\end{proof}

The following proposition was proved by White \cite{White2001} and
Lam~ \cite{Lam2004} for $\mu=\emptyset$ and $\mu=(1)$. In our
definition of $\s(D)$, it holds for any $\mu$. Our proof is similar to
Lam's~ \cite[Proposition 21]{Lam2004}. Recall that a SDT is a SYT with
the condition that two cells with entries $2i-1$ and $2i$ make a
domino. Thus the sign of a SDT is just the sign of a SYT.

\begin{prop}\Label{ev}
Let $D$ be a SDT of shape $\lm$. Then 
$$\s(D)=(-1)^{ev(D)}.$$
\end{prop}
\begin{proof}
We use induction on $n$, the number of dominoes in $D$. It is trivial
if $n=0$. Let $sh(D)=\lm\vdash 2n$. Let ${\bf d}$ be the domino with
entry $n$ and let $a$ and $b$ be the cells in ${\bf d}$ with $a\lhd
b$.  Let $D'$ be the SDT obtained from $D$ by removing ${\bf d}$. Let
$T_0\in\t(\mu)$.  Then $\s(D) = \s(T_0)(-1)^{\inv(T_0\diamond D)}$ and
$\s(D') = \s(T_0)(-1)^{\inv(T_0\diamond D')}$.

Since $(T_0\diamond D)(a)$ and $(T_0\diamond D)(b)$ are greater than
any entry of $T_0\diamond D'$,
$$\Inv(T_0\diamond D)=\Inv(T_0\diamond D') \cup\{(a,c) : a\lhd c,
c\in\lambda\setminus {\bf d}\} \cup\{(b,c) : b\lhd c, c\in\lambda
\setminus {\bf d}\}.$$ Thus we have
$$\inv(T_0\diamond D)\equiv \inv(T_0\diamond D')+\#\{c\in \lambda:a\lhd
c\lhd b\} \mod 2.$$ If ${\bf d}$ is horizontal then $\#\{c\in
\lambda:a\lhd c\lhd b\}=0$. If ${\bf d}$ is vertical in the $i$-th
column then $\#\{c\in \lambda:a\lhd c\lhd b\}=i-1$. Thus
$$\#\{c\in \lambda:a\lhd c\lhd b\} \equiv ev(D)-ev(D') \mod 2.$$ Since
$\s(D')=(-1)^{ev(D')}$ by the induction hypothesis, we get
\begin{align*}
\s(D)&=\s(T_0)(-1)^{\inv(T_0\diamond D)} =
\s(T_0)(-1)^{\inv(T_0\diamond D')+ev(D)-ev(D')}\\ &=
\s(D')(-1)^{ev(D)-ev(D')}=(-1)^{ev(D)}. \qedhere
\end{align*}
\end{proof}

\subsection{Sign-imbalance of skew shapes}
The {\em \si.} $I_{\lm}$ of a \ss. $\lm$ is defined by
$$I_{\lm}=\sum_{T\in\t(\lm)} \s(T).$$ 

Let $\lm \vdash 2n$ and $T\in\t(\lm)$. If $2k-1$ and $2k$ are neither
in the same row nor in the same column of $T$ for some $k$, let $T'$
be the \SYT. obtained from $T$ by switching the entries $2k-1$ and
$2k$ for the smallest such $k$.  Then $T\mapsto T'$ is a sign
reversing \iv. on $\t(\lm)\setminus\d(\lm)$. Thus we only need to
consider SDTs.  Then, using \autoref{ev}, we get
$$I_{\lm}=\sum_{D\in \d(\lm)}\s(D)=\sum_{D\in \d(\lm)}(-1)^{ev(D)}.$$

The idea of the following lemma is found in the proof of Corollary~24
in Lam's paper \cite{Lam2004}.
\begin{lem}\Label{iandf}
  Let $n\geq0$ and $\lm\vdash 2n$. Then
$$I_{\lm} = (-1)^{-\frac12 (\frac12|\lm| - h(\lm))} f_2^{\lm}(-1).$$
\end{lem}
\begin{proof}
Using the above argument and \autoref{stat},
\begin{align*}
I_{\lm} &= \sum_{D\in\d(\lm)}(-1)^{ev(D)}
=\sum_{D\in\d(\lm)}(-1)^{-\frac12(ov(D)-ev(D)) +sp(D)}\\
&=(-1)^{-\frac12 (\frac12|\lm| - h(\lm))} f_2^{\lm}(-1). \qedhere
\end{align*}
\end{proof}

Note that, in the above lemma, although both $(-1)^{-\frac12
  (\frac12|\lm| - h(\lm))}$ and $f_2^{\lm}(-1)$ lie in
$\Z[\sqrt{-1}]$, it is easy to check that their product is an integer.

Now we get a generalization of Eq.~\eqref{two} to \ss.s of even
size. In \autoref{sjthm}, we prove a stronger theorem which has no
restriction on the size of \ss.s.
\begin{cor} 
Let $\alpha$ and $\beta$ be fixed partitions and $n$ and $m$ be fixed
nonnegative integers. Then
$$\sum_{\substack{\la\vdash 2m \\ \lb\vdash 2n}} (-1)^{v(\lambda)}
I_{\la}I_{\lb} = (-1)^{v(\alpha)+v(\beta)} \sum_{ \substack{\bm\vdash
    2m \\ \am \vdash 2n}} (-1)^{v(\mu)} I_{\bm} I_{\am}.$$
\end{cor}
\begin{proof}
If $q=-1$ in \autoref{double}, then
$$\sum_{\substack{\la\vdash 2m \\ \lb\vdash 2n}} f^{\la}_2(-1) f^{\lb}_2(-1)
  = \sum_{ \substack{ \bm\vdash 2m \\ \am \vdash 2n}} f^{\bm}_2(-1)
  f^{\am}_2(-1).$$
Let $\eta(\lambda) = \frac12|\lambda/\tl|
  -h(\lambda/\tl)$.  Then for a \ss. $\lm$ with $\tl=\tm$ we have
$$\frac{|\lm|}2 - h(\lm) = \eta(\lambda) - \eta(\mu).$$ Since we can
  assume $\tl = \tm = \tilde{\alpha} = \tilde{\beta}$ (or
  equivalently, $\la$, $\lb$, $\bm$ and $\am$ are \dtl.), by
  \autoref{iandf} we get
$$\sum_{\substack{\la\vdash 2m \\ \lb\vdash 2n}} (-1)^{\eta(\lambda)}
I_{\la}I_{\lb} = (-1)^{\eta(\alpha)+\eta(\beta)} \sum_{\substack{
    \bm\vdash 2m \\ \am \vdash 2n}} (-1)^{\eta(\mu)} I_{\bm}I_{\am}.$$
By \autoref{stat}, we have $\eta(\lambda) \equiv v(\lambda/\tl) \mod
2$, which finishes the proof.
\end{proof}

\subsection{Definition of a generalized \si. formula}
Let $\alpha$ be a fixed partition and $n\geq0$.
We define
$$F_n^{\alpha} = F_n^{\alpha}(x,y,z)=\sum_{\la\vdash n}
x^{v(\la)}y^{h(\la)}z^{d(\la)} I_{\la}.$$ Then Eq.~\eqref{one} can be
written as $F_{n}^{\emptyset}(x,y,z)=(x+y)^{\flr n}$.

Let $\alpha^+$ denote the set $\{\lambda: |\lambda|=|\alpha|+1,
\alpha\subset \lambda\}$. For $\lambda\in\alpha^+$, let
$u(\lambda,\alpha)$ denote the number of cells $a\in\alpha$ such that
$b\lhd a$ for the unique cell $b\in \la$. For example, if
$\alpha=(7,5,5,2)$ and $\lambda=(7,6,5,2)$ then
$u(\lambda,\alpha)=7$.

\def\na{\nu/\alpha}

\begin{prop}\Label{decrease}
Let $\alpha$ be a fixed partition and $n\geq0$. Then
$$F_{n+1}^{\alpha}= \sum_{\nu\in \alpha^+}(-1)^{u(\nu,\alpha)} \pp{\na}
F_{n}^{\nu},$$ where $\pp{\na} = x^{v(\na)}y^{h(\na)}z^{d(\na)}$.
\end{prop}
\begin{proof}
Let $\la\vdash n+1$.  If $T\in\t(\la)$ then the cell whose entry is 1
must be the unique cell of $\na$ for some $\nu\in \alpha^+$. Since
$\na$ contains only one cell, there is a unique \SYT. of shape $\na$,
say $T_{\nu}$. Then $\s(T_{\nu})=(-1)^{u(\nu,\alpha)}$. Thus
$T\in\t(\la)$ \iff. $T=T_{\nu}\diamond T'$ for some $\nu\in \alpha^+$
and $T'\in\t(\ln)$, which implies
$$I_{\la} = \sum_{\nu\in\alpha^+}
\sum_{T'\in\t(\ln)}\s(T_{\nu}\diamond T') = \sum_{\nu\in\alpha^+}
(-1)^{u(\nu,\alpha)} I_{\ln}.$$ Since
$\pp{\la}=\pp{\na}\cdot\pp{\ln}$,
\begin{align*}
F_{n+1}^{\alpha} &= \sum_{\la\vdash n+1} \pp{\la}
\sum_{\nu\in\alpha^+} (-1)^{u(\nu,\alpha)} I_{\ln}\\ &=
\sum_{\nu\in\alpha^+} (-1)^{u(\nu,\alpha)} \pp{\na} \sum_{\ln\vdash n}
\pp{\ln} I_{\ln}\\ &= \sum_{\nu\in \alpha^+}(-1)^{u(\nu,\alpha)}
\pp{\na} F_{n}^{\nu}. \qedhere
\end{align*}
\end{proof}
Using \autoref{decrease}, we can calculate $F_{n}^{\alpha}$ for all
$n\geq0$ if we have $F_{n}^{\alpha}$ for all even $n$.  Thus we will
focus on skew shapes $\lm\vdash 2n$ of even size.

We extend the definition of the \si. $I_{\lm}$ to reversed shapes as
follows.  For a reversed shape $\lm\vdash -2n$, define
$$I_{\lm} = (-1)^{-\frac12 (\frac12|\lm| - h(\lm))} f_2^{\lm}(-1).$$
Note that the above equation is the same one in \autoref{iandf}. We
have a relation between $I_{\lm}$ and $I_{\ml}$.
\begin{prop}\Label{fi}
  Let $\lm\vdash -2n$ be a \rs. for $n\geq0$. Then
$$I_{\lm} = (-1)^{v(\ml)}I_{\ml}.$$
\end{prop}
\begin{proof} 
If $\ml$ is not \dtl. then $I_{\lm} = I_{\ml}=0$.
Otherwise, we have $n-h(\ml)\equiv v(\ml) \mod
2$, by \autoref{stat}. Thus,
\begin{align*}
I_{\lm} &= (-1)^{-\frac12(-n-h(\lm))} f_2^{\lm}(-1)\\
 &= (-1)^{n-h(\ml)}(-1)^{-\frac12(n- h(\ml))} f_2^{\ml}(-1)\\
 &= (-1)^{n-h(\ml)}I_{\ml} = (-1)^{v(\ml)}I_{\ml}. \qedhere
\end{align*}
\end{proof}

Now we extend the definition of $F_{2n}^{\alpha}$ as follows: for
$n\geq0$, define
$$F_{-2n}^{\alpha} = \sum_{\la\vdash -2n}
x^{v(\la)}y^{h(\la)}z^{d(\la)}I_{\la}.$$ Then, by \autoref{fi},
$$F_{-2n}^{\alpha} = \sum_{\al\vdash 2n}
(-x)^{-v(\al)}y^{-h(\al)}z^{-d(\al)}I_{\al}.$$ 

\subsection{A method to obtain a generalized \si. formula}

\begin{lem}
  Let $\alpha$ be a fixed partition and $n\geq0$. Then
$$F_{2n}^{\alpha} = W_n^{\alpha} \left( (x\sqrt z)\inverse, (y\sqrt
z\sqrt{-1})\inverse, -1 \right) \cdot (xy\sqrt z)^n,$$  
and
$$F_{-2n}^{\alpha} = W_{-n}^{\alpha} \left( (x\sqrt z)\inverse, (y\sqrt
z\sqrt{-1})\inverse, -1 \right) \cdot (xy\sqrt z)^{-n}.$$  
\end{lem}
\begin{proof}
  Let $\la\vdash 2n$ be a \dtl. \ss.. Then, by \autoref{stat}, we have
  $\frac12(v(\la)+h(\la)-\frac12|\la|) = d(\la)$. Thus
  \begin{align*}
   & W_n^{\alpha} \left( (x\sqrt z)\inverse, (y\sqrt z\sqrt{-1})\inverse,
    -1 \right) \cdot (xy\sqrt z)^n\\ &= \sum_{\la\vdash 2n} (x\sqrt
    z)^{v(\la)-\frac12|\la|} (y\sqrt z\sqrt{-1})^{h(\la)-\frac12|\la|}
    f_2^{\la}(-1)\cdot(xy\sqrt z)^{\frac12|\la|} \\ &= \sum_{\la\vdash 2n}
    x^{v(\la)} y^{h(\la)} z^{\frac12(v(\la)+h(\la)-\frac12|\la|)}
      (-1)^{\frac12(h(\la)-\frac12|\la|)} f_2^{\la}(-1)\\ &= \sum_{\la\vdash 2n}
      x^{v(\la)} y^{h(\la)} z^{d(\la)} I_{\la} = F_{2n}^{\alpha}.
  \end{align*}

  Now let $\la\vdash -2n$ be a reversed shape such that $\al$ is
  domino-tileable.  Since all the arguments we used here are remaining
  true when we change $n$ to $-n$, we get the second identity in the
  lemma as well.
\end{proof}

Now we get a \gf. for $F_{2n}^{\alpha}$.
\begin{thm}\Label{main} 
 Let $\alpha$ be a fixed partition. Then
$$\frac{\sum_{n\geq0}F_{2n}^{\alpha} \frac{t^n}{n!}}  {\sum_{n\geq0}
    F_{-2n}^{\alpha} \frac{(x^2y^2zt)^n}{n!}}  = \exp\left( (x+y)t
  \right).$$
\end{thm}
\begin{proof}
Substitute $x,y,q$ and $t$ in \autoref{gf} with $(x\sqrt z)\inverse,
(y\sqrt z\sqrt{-1})\inverse, -1$ and $xy\sqrt z t$.  Then we get this
theorem.
\end{proof}

\begin{cor}\Label{F2n}
  Let $\alpha$ be a fixed partition with $\alpha/\ta\vdash 2k$. Then,
  for $n\geq0$,
$$F_{2n}^{\alpha}=\sum_{j=0}^k \binom{n}{j} (x+y)^{n-j} (x^2y^2z)^j
F_{-2j}^{\alpha}.$$
\end{cor}
If $\alpha=\dr$ then $\alpha/\ta\vdash 0$. Thus we get the following
corollary.
\begin{cor}\Label{xyn}
For any integers $k\geq0$ and $n\geq0$, we have $F_{2n}^{\dk}=(x+y)^n$.
\end{cor}
The next example shows how to calculate $F_{2n}^{\alpha}$.
\begin{example}
Let us find $F_{2n}^{\alpha}$ for $\alpha=(2,2)$.  We have
$\ta=\emptyset$ and $\alpha/\ta\vdash 4$. Using Table~\ref{tab:st} we
get
\begin{align*}
  F_{0}^{(2,2)} &= 1, \\
  F_{-2}^{(2,2)} &= (-x)^{-2}y^{-1}z^{-1} + (-x)^{-1}y^{-2}z^{-1}(-1)
=x^{-2}y^{-2}z^{-1}(x+y), \\
  F_{-4}^{(2,2)} &= 0.
\end{align*}

\begin{table}
  \begin{center}
\begin{tabular}{|c|c|c|c|c|} \hline
  $j$ & 0 & \multicolumn{2}{c|}{1} & 2 \\ \hline
  $\lambda$ & $(2,2)$ & $(2)$ & $(1,1)$ & $\emptyset$ \\ \hline
  $\al$
\psraise (5,2){\pspicture (0,-2.5) (0,0) \endpspicture}
 & 
\psraise (2,1){\pspicture (0,-2) (2,0)\pspolygon[linestyle=dotted] (2,0) (2,-2) (0,-2) (0,0) \endpspicture}
 & 
\psraise (2,1){\pspicture (0,-2) (2,0) \cell(2,1)[] \cell(2,2)[]\pspolygon[linestyle=dotted] (2,0) (2,-1) (0,-1) (0,0) \endpspicture}
           & 
\psraise (2,1){\pspicture (0,-2) (2,0) \cell(1,2)[] \cell(2,2)[]\pspolygon[linestyle=dotted] (1,0) (1,-2) (0,-2) (0,0) \endpspicture}
 & 
\psraise (2,1){\pspicture (0,-2) (2,0) \cell(1,1)[] \cell(1,2)[] \cell(2,1)[] \cell(2,2)[] \endpspicture}
 \\ \hline
  $v(\al)$ & $0$ & $2$ & $1$ & $2$ \\ \hline
  $h(\al)$ & $0$ & $1$ & $2$ & $2$ \\ \hline
  $d(\al)$ & $0$ & $1$ & $1$ & $1$ \\ \hline
  $I_{\al}$ & $1$ & $1$ & $-1$ & $0$ \\ \hline
\end{tabular}
  \end{center}
\caption{Statistics of $\al$ for $\alpha=(2,2)$,
  $\lambda$ and $j$ with $\al\vdash 2j$.}\label{tab:st}
\end{table}

Thus $F_{2n}^{(2,2)} = (n+1)(x+y)^n$ and
$$\sum_{n\geq0}F_{2n}^{(2,2)} \frac{t^n}{n!} = \left( 1 + (x+y)t
\right) \cdot \exp\left( (x+y)t \right).$$
\end{example}

\subsection{A closed formula for a staircase partition}
Now we can get a closed formula for $F_{n}^{\dk}$.
\begin{thm}\Label{fnk}
For any integers $k\geq0$ and $n\geq0$, we have
\begin{align*}
F_{2n}^{\dk} &= (x+y)^{n},\\ F_{2n+1}^{\dk} &= \left
\{ \begin{array}{ll} (x+y)^{n}, & \mbox{if } k\equiv 0 \mod 4,
  \\ (x+y)^{n+1}, & \mbox{if } k\equiv 1 \mod 4, \\ xyz(x+y)^{n}, &
  \mbox{if } k\equiv 2 \mod 4, \\ 0, &\mbox{if } k\equiv 3 \mod 4.\\
  \end{array}\right.
\end{align*}
\end{thm}
We have already proved the even case in \autoref{F2n}. For the odd
case we need two lemmas.  For $0\leq i\leq k$, let $\dk^i$ denote the
partition in $\dk^+$ obtained from $\dk$ by adding the cell
$(k+1-i,i+1)$.  Recall $\pp{\lm} = x^{v(\lm)}y^{h(\lm)}z^{d(\lm)}$,
which is used in \autoref{decrease}.

\begin{lem}\Label{lem2}
Let $k\geq0$. Then
  $$ \sum_{i=0}^k (-1)^{\flr i} \pp{\dk^i/\dk} = \frac{1+(-1)^{\flr
      k}}2 x^{\frac{1-(-1)^k}2} + \frac{1+(-1)^{\flr {k-1}}}2
  x^{\frac{1+(-1)^k}2}y z^{\frac{1+(-1)^k}2}.$$
\end{lem}
\begin{proof}
Since $\dk^i/\dk$ contains only one cell $(k-i+1,i+1)$, we have
\begin{align*}
\sum_{i=0}^k (-1)^{\flr i} \pp{\dk^i/\dk} &= \sum_{i=0}^k
(-1)^{\flr i} x^{\frac{1+(-1)^{k-i+1}}2} y^{\frac{1+(-1)^{i+1}}2} 
z^{\frac{1+(-1)^{k-i+1}}2 \cdot \frac{1+(-1)^{i+1}}2}\\ &= \sum_{a=0}^{\flr
  k} (-1)^a x^{\frac{1-(-1)^{k}}2} + \sum_{b=0}^{\flr{k-1}} 
(-1)^b x^{\frac{1+(-1)^{k}}2}yz^{\frac{1+(-1)^{k}}2}.
\end{align*}
 Since $\sum_{i=0}^m (-1)^i= \frac{1+(-1)^{m}}2$, we are done.
\end{proof}

\begin{lem}\Label{lem3}
Let $j\geq1$ and $\lambda$ be a fixed partition with
$\dk^i/\lambda\vdash 2j$ for some $i=0,1,\ldots,k$. Then
$$\sum_{i=0}^k (-1)^{\flr i + i} I_{\dk^i/\lambda}=0.$$
\end{lem}
\begin{proof}
\def\dd{\mathfrak{D}} Let $\dd=\cup_{i=0}^k \d(\dk^i/\lambda)$. For
$D\in\dd$ with $sh(D)=\dk^i/\lambda$, we define $s(D)=(-1)^{\flr i +
  i+ev(D)}$. Since
$$\sum_{i=0}^k (-1)^{\flr i + i} I_{\dk^i/\lambda} = \sum_{D\in\dd}
s(D),$$ it is sufficient to construct an involution
$\omega:\dd\rightarrow\dd$ satisfying $s(\omega(D))=-s(D)$.  Let
$D\in\dd$ and ${\bf d}$ be the domino of $D$ with the largest entry.
Then ${\bf d}\cap \dk$ must have only one cell, say $(a,b)$.  Let
${\bf d}'$ be the domino satisfying ${\bf d}\cup {\bf d}'
=\{(a,b),(a+1,b),(a,b+1)\}$. We define $\omega(D)$ to be the SDT
obtained from $D$ by changing ${\bf d}$ with ${\bf d}'$,
see Fig.~\ref{omega}. It is
obvious that $\omega$ is an involution. 

\begin{figure}
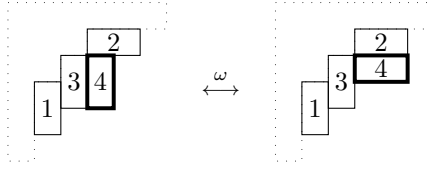

  \begin{center}
\psraise (6,1){\pspicture (0,-6) (6,0)\pspolygon[linestyle=dotted] (6,0) (6,-1) (3,-1) (3,-2) (2,-2) (2,-3) (1,-3) (1,-6) (0,-6) (0,0)\pspolygon (2,-3) (2,-5) (1,-5) (1,-3)\rput (1.5,-4){1}\pspolygon (3,-2) (5,-2) (5,-1) (3,-1)\rput (4,-1.5){2}\pspolygon (2,-4) (3,-4) (3,-2) (2,-2)\rput (2.5,-3){3}\pspolygon (3,-4) (4,-4) (4,-2) (3,-2)\rput (3.5,-3){4}\pspolygon[linewidth=1.5pt] (3,-4) (4,-4) (4,-2) (3,-2) \endpspicture}
 \quad $\substack{\omega\\ \longleftrightarrow}$
 \quad 
\psraise (6,1){\pspicture (0,-6) (6,0)\pspolygon[linestyle=dotted] (6,0) (6,-1) (3,-1) (3,-2) (2,-2) (2,-3) (1,-3) (1,-6) (0,-6) (0,0)\pspolygon (2,-3) (2,-5) (1,-5) (1,-3)\rput (1.5,-4){1}\pspolygon (3,-2) (5,-2) (5,-1) (3,-1)\rput (4,-1.5){2}\pspolygon (2,-4) (3,-4) (3,-2) (2,-2)\rput (2.5,-3){3}\pspolygon (3,-3) (5,-3) (5,-2) (3,-2)\rput (4,-2.5){4}\pspolygon[linewidth=1.5pt] (3,-3) (5,-3) (5,-2) (3,-2) \endpspicture}
  \end{center}
  \caption{Involution $\omega$.}\label{omega}
\end{figure}

To show $s(\omega(D))=-s(D)$, we can assume that $sh(D)=\dk^i/\lambda$
and ${\bf d}$ is a vertical domino.  Then $sh(\omega(D))=
\dk^{i+1}/\lambda$ and $ev(D)-ev(\omega(D))$ is 1 if $i$ is odd, and
$0$ if $i$ is even, which implies
$(-1)^{ev(D)-ev(\omega(D))}=(-1)^{i}$. Then we get
$$\frac{s(D)}{s(\omega(D))}=(-1)^{\flr i + i+ev(D) - (\flr {i+1} +
  i+1+ev(\omega(D)))} =(-1)^{\flr i - \flr {i+1} + i-1}=-1,$$ 
which completes the proof.
\end{proof}

\begin{proof}[Proof of \autoref{fnk}]
We will show the following equivalent equation:
$$F_{2n+1}^{\dk} = (x+y)^n \left( \frac{1+(-1)^{\flr k}}2
x^{\frac{1-(-1)^k}2} + \frac{1+(-1)^{\flr {k-1}}}2
x^{\frac{1+(-1)^k}2}y z^{\frac{1+(-1)^k}2} \right).$$ Since
$\dk^+=\{\dk^i : i=0,1,2,\ldots,k\}$, by \autoref{decrease}
and \autoref{F2n},
\begin{align*}
F_{2n+1}^{\dk} &= \sum_{i=0}^k (-1)^{u(\dk^i/\dk)} \pp{\dk^i/\dk}F_{2n}^{\dk^i}\\
&= 
\sum_{j=0}^k \binom{n}{j} (x+y)^{n-j} (x^2y^2z)^j
\sum_{i=0}^k (-1)^{\flr i} \pp{\dk^i/\dk} F_{-2j}^{\dk^i}.
\end{align*}
By \autoref{lem2}, it is sufficient to show that, for $j\geq 1$, the
following sum is 0:
\begin{align*}
&\sum_{i=0}^k (-1)^{\flr i} \pp{\dk^i/\dk}F_{-2j}^{\dk^i} \\
&= \sum_{i=0}^k (-1)^{\flr i} \pp{\dk^i/\dk}\sum_{\dk^i/\lambda\vdash 2j}
(-x)^{-v(\dk^i/\lambda)} y^{-h(\dk^i/\lambda)} z^{-d(\dk^i/\lambda)} I_{\dk^i/\lambda}\\
&= \sum_{i=0}^k (-1)^{\flr i} 
\sum_{\dk^i/\lambda\vdash 2j} (-x)^{-v(\dk/\lambda)} y^{-h(\dk/\lambda)}
z^{-d(\dk/\lambda)} (-1)^{v(\dk^i/\dk)} I_{\dk^i/\lambda}\\
&= \sum_{i=0}^k (-1)^{\flr i} 
\sum_{\dk^i/\lambda\vdash 2j} (-x)^{-v(\dk/\lambda)} y^{-h(\dk/\lambda)}
z^{-d(\dk/\lambda)} (-1)^{k-i} I_{\dk^i/\lambda}\\
&= \sum_{\lambda}(-x)^{-v(\dk/\lambda)} y^{-h(\dk/\lambda)}
z^{-d(\dk/\lambda)} (-1)^{k}
\sum_{i=0}^k (-1)^{\flr i + i} I_{\dk^i/\lambda},
\end{align*}
where the last sum is over $\{\lambda:\dk^i/\lambda\vdash 2j$
for some $i\}$. By \autoref{lem3}, we are done.
\end{proof}

\section{Generalizing \sj.'s theorems}\label{sjthm}

In this section, we only consider (skew) SYTs.

Let $\m$ denote the set of $n\times m$ matrices whose entries are $0$
or $1$ such that the total number of $1$'s is $j$ and there is at most
one $1$ in each row and column.  For $M\in\m$, let $perm(M)$ denote
the permutation whose permutation matrix is obtained from $M$ by
removing the rows and columns without $1$'s.

Using the same argument of \autoref{gdi} with the usual local rules
for the \RS. correspondence, we can formulate the following theorem.

\begin{thm}\Label{gdih} \cite{Roby1991}
Let $\alpha$ and $\beta$ be fixed partitions and $n$ and $m$ be fixed
nonnegative integers.  Then $\gdi=\up\circ(\down)\inverse$ induces a
bijection
$$\gdi : \bigcup_{j\geq0} \left( \bigcup_{\substack{\bm\vdash m-j
    \\ \am\vdash n-j}} \t(\bm)\times\t(\am)\times \m
\right) \rightarrow \bigcup_{\substack{\la\vdash m \\ \lb\vdash n}}
\t(\la)\times\t(\lb).$$
\end{thm}

The following elegant theorem was proved by Reifegerste
\cite{Reifegerste2004} and \sj. \cite{Sjostrand2005} independently.
\begin{thm}\Label{reif}
Let $\pi$ correspond to $(P,Q)$ in the \RS. correspondence and
$sh(P)=\lambda$. Then
$$\s(\pi)=(-1)^{v(\lambda)}\s(P)\s(Q).$$
\end{thm}

By the local rules, the next lemma is an immediate result of
\autoref{reif}.
\begin{lem}\Label{emp}
Let $P,Q\in\t(\lambda)$ and $M\in\m$ satisfy
$\gdi(\emptyset_{\emptyset},\emptyset_{\emptyset},M)=(P,Q)$, where
$\emptyset_{\emptyset}$ denotes the empty \SYT. of shape
$\emptyset/\emptyset$.  Then
$$\s(perm(M))=(-1)^{v(\lambda)}\s(P)\s(Q).$$
\end{lem}

Let $k=n+m-j$.  For $M\in\m$, let $\overline{M}$ denote the element in
$\mathfrak{m}_{k,k}^{k}$ which can be expressed as
$\left(\substack{A\\C} \substack{B\\M}\right)$ such that $A={\bf 0}$,
$perm(B)=12\cdots (m-j)$ and $perm(C) =12\cdots (n-j)$. It is easy to
check that such $\overline{M}$ exists uniquely.  For example, if
$M=\mat{cccc}{ 0&1&0&0\\ 0&0&0&0\\0&0&0&0}$ then
$\overline{M}=\mat{cccccc}{
  0&0&1&0&0&0\\
  0&0&0&0&1&0\\
  0&0&0&0&0&1\\
  0&0&0&1&0&0\\ 1&0&0&0&0&0\\ 0&1&0&0&0&0}$.

For a permutation $\pi$, let $\inv(\pi)$ denote the number of
inversions, i.e., pairs $(i,j)$ such that $i<j$ and $\pi_i>\pi_j$.
Let $\inv(M)=\inv(perm(\overline{M}))$.  The sign of $M$ is defined by
$$\s(M)=(-1)^{\inv(M)}.$$ 

For nonnegative integers $n$, $k$ and $\seq{a}{r}{,}$ such that
$\sum_{i=1}^r a_i=n$, we denote
$$\qfactorial{n} =(1+q)(1+q+q^2)\cdots(1+q+\cdots+q^{n-1}),$$

$$\qbinom{n}{\seq{a}{r}{,}}=
\frac{\qfactorial{n}}{\qfactorial{a_1}\qfactorial{a_2}
  \cdots\qfactorial{a_r}},\quad \qbinom{n}{k} = \qbinom{n}{k, n-k}.$$
\begin{prop}\Label{invM}
Let $n$, $m$ and $j$ be nonnegative integers. Then
  $$\sum_{M\in\m}q^{\inv(M)}=q^{(n-j)(m-j)}
  \qbinom{n}{j}\qbinom{m}{j}\qfactorial{j}.$$
\end{prop}
\begin{proof}
Let $M\in\m$ and $\pi=perm(M)$. Let ${\bf r}=\seq{r}{n}{}$
(resp. ${\bf c}=\seq{c}{m}{}$) be the $(0,1)$-sequence such that
$r_i=0$ (resp. $c_i=0$) \iff. the $i$-th row (resp. column) of $M$
contains 1. Then ${\bf r}\in\S(\{0^{j},1^{n-j}\})$, ${\bf
  c}\in\S(\{0^{j},1^{m-j}\})$ and $\pi\in\S([j])$, where $\S(X)$
denotes the set of permutations of $X$ for a (multi)set $X$.  It is
well known, for example see \cite[Proposition 1.3.17]{Stanley1997},
that if $X=\{1^{a_1},2^{a_2},\ldots,n^{a_n}\}$ then
$\sum_{\pi\in\S(X)}q^{\inv(\pi)}=\qbinom{n}{\seq{a}{n}{,}}$.  Since
$\inv(M)=(n-j)(m-j)+\inv({\bf r})+\inv({\bf c})+\inv(\pi)$, we are
done.
\end{proof}

Substituting $q=-1$ in \autoref{invM}, we get the following lemma.
\begin{lem}\Label{signM}
Let $n$, $m$ and $j$ be nonnegative integers. Then
  $$\sum_{M\in\m}\s(M)=\left\{
  \begin{array}{ll}
    (-1)^{mn}, & \mbox{if $j=0$,}\\
    \frac{1-(-1)^{mn}}2, & \mbox{if $j=1$,}\\
    0, & \mbox{otherwise.}
  \end{array}
  \right.$$
\end{lem}

Now we can generalize Eq.~\eqref{sjeq}.
\begin{thm}\Label{sjsign}
Let $U\in\t(\bm)$, $V\in\t(\am)$, $P\in\t(\la)$, $Q\in\t(\lb)$ and a
matrix $M$ satisfy $\gdi(U,V,M)=(P,Q)$.  Then
$$(-1)^{v(\alpha)+v(\beta)+v(\lambda)}\s(P)\s(Q)=(-1)^{v(\mu)}
\s(U)\s(V)\s(M).$$
\end{thm}
\begin{proof}
Let $\la\vdash m$, $\lb\vdash n$ and $\lambda\vdash k$.  Let
$A\in\t(\alpha)$ and $B\in\t(\beta)$. Then there is a unique $k\times
k$ full \gd. $G=(\Gamma,N)$ with $\up(G)= (A\diamond P, B\diamond
Q)$. It is obvious that $N=\left(\substack{M_{11}\\M_{21}}
\substack{M_{12}\\M}\right)$ for suitable matrices $M_{11}$, $M_{12}$
and $M_{21}$.  We can construct \gd.s $G_{11}$, $G_{21}$ and
$G_{12}$ from $G$ as follows:
\begin{align*}
G_{11} &= \left( \left(\Gamma_{(i,j)}\right)_{
\begin{subarray}{l}
    0\leq i\leq k-n\\0\leq j\leq k-m\end{subarray}
},
 M_{11} \right),\\
G_{21} &= \left(
\left(\Gamma_{(i,j)}\right)_{
  \begin{subarray}{l}
    0\leq i\leq k\\0\leq j\leq k-m
  \end{subarray}
}, \left(\substack{M_{11}\\M_{21}}\right) \right),\\ 
G_{12} &=
\left(\left( \Gamma_{(i,j)}\right)_{
  \begin{subarray}{l}
    0\leq i\leq k-n\\0\leq j\leq k
  \end{subarray}
}, (M_{11} M_{12})\right).
\end{align*}
Let $U_0=(G_{11})\bot ^{\rm SYT}\in\t(\mu)$ and $V_0=(G_{11})\rig
^{\rm SYT}\in\t(\mu)$, where $C^{\rm SYT}$ is defined similarly to
$C^{\rm SDT}$ in \autoref{growth}.  See Fig.~\ref{fig:ex}, which
roughly represents $G$ and these \SYT.s.
\begin{figure}
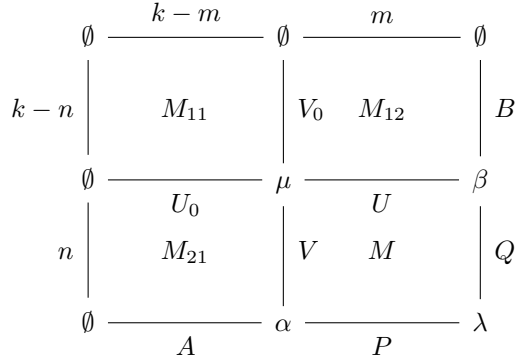

  \begin{center}
$\begin{psmatrix}[nodesep=0.5,rowsep=1.5,colsep=2.5]
\emptyset &       & \emptyset &       &\emptyset \\
          & M_{11} &           & M_{12} &          \\
\emptyset &       &    \mu    &       &\beta     \\
          & M_{21} &           & M     &          \\
\emptyset &       &    \alpha &       &\lambda
\ncline{1,1}{1,3}^{k-m} \ncline{1,3}{1,5}^{m} \ncline{1,1}{3,1}<{k-n}
\ncline{3,1}{5,1}<{n} \ncline{1,3}{3,3}>{V_0} \ncline{1,5}{3,5}>{B}
\ncline{3,3}{5,3}>{V} \ncline{3,5}{5,5}>{Q} \ncline{3,1}{3,3}_{U_0}
\ncline{3,3}{3,5}_{U} \ncline{5,1}{5,3}_{A} \ncline{5,3}{5,5}_{P}
\end{psmatrix}$
  \end{center}
\caption{Growth diagrams and \SYT.s.}\label{fig:ex}
\end{figure}

Let $perm(M_{11})=\gamma$,
$cp\left(\substack{M_{11}\\M_{21}}\right)=\sigma$, $perm(M_{11}
M_{12})=\tau$ and $cp\left(\substack{M_{11}\\M_{21}}
\substack{M_{12}\\M} \right)=\pi$. Then $\s(\pi) =
\s(\sigma)\s(\tau)\s(\gamma)\s(M)$, and by \autoref{emp},
\begin{align*}
\s(\pi) &= (-1)^{v(\lambda)}\s(A\diamond P)\s(B\diamond Q),\\
\s(\sigma) &= (-1)^{v(\alpha)}\s(A)\s(V_0\diamond V),\\
\s(\tau) &= (-1)^{v(\beta)}\s(U_0\diamond U)\s(B),\\
\s(\gamma) &= (-1)^{v(\mu)}\s(U_0)\s(V_0).
\end{align*}
Multiplying the above five equations, we get this theorem.
\end{proof}
\begin{remark}
\sj.'s theorem, which is Eq.~\eqref{sjeq}, is stated in a different
way, however, it is not difficult to see that it is equivalent to
\autoref{sjsign} with $\alpha=\beta$. Also note that, \sj. used $\s_1$
for the sign of a \SYT..  Despite the different definitions, by
\autoref{sign}, if $sh(P)=sh(Q)$ then $\s_1(P)\s_1(Q)=\s(P)\s(Q)$.
\end{remark}

Using \autoref{gdih}, \autoref{signM} and \autoref{sjsign}, we get the
following generalization of Eq.~\eqref{sjsi}.
\begin{thm}\Label{sjss}
Let $\alpha$ and $\beta$ be fixed partitions and $n$ and $m$ be fixed
nonnegative integers. Then
  \begin{align*}
 &(-1)^{v(\alpha)+v(\beta)} \sum_{\substack{\la\vdash m \\ \lb\vdash
        n}} (-1)^{v(\lambda)} I_{\la}I_{\lb}\\ &= (-1)^{mn}\sum_{
      \substack{\bm\vdash m \\ \am \vdash n}} (-1)^{v(\mu)} I_{\bm}
    I_{\am} +\frac{1-(-1)^{mn}}2 \sum_{ \substack{\bm\vdash m-1 \\ \am
        \vdash n-1}} (-1)^{v(\mu)} I_{\bm} I_{\am}.
  \end{align*}
\end{thm}

\section*{Acknowledgments}
I would like to thank my advisor, Dongsu Kim, for his encouragement
and helpful comments. I would also like to thank Seunghyun Seo for his
careful reading of the first version of this paper.

\end{document}